\documentclass [a4,11pt]{amsart}

\usepackage{amsmath,amssymb,amsfonts,enumerate,amsthm, amscd,}

\newtheorem{thm}{Theorem}[section]
\newtheorem{cor}[thm]{Corollary}
\newtheorem{lem}[thm]{Lemma}

\newtheorem{exam}[thm]{Example}
\newtheorem{rem}[thm]{Remark}

\def\proof{{\parindent0pt {\bf Proof.\ }}}
\newcommand{\field}[1]{\mathbb{#1}}

\newcommand{\Q }{\field{Q}}
\newcommand{\Z }{\field{Z}}

\theoremstyle{definition}

\theoremstyle{remark}

\theoremstyle{Definition and Notation}

\begin{document}

\bibliographystyle{amsplain}


\title[Coherence in amalgamated algebra along an ideal]{Coherence in amalgamated algebra along an ideal}

\author{Karima Alaoui Ismaili}
\address{Karima Alaoui Ismaili\\Department of Mathematics, Faculty of Science and Technology of Fez, Box 2202,
University S.M. Ben Abdellah Fez, Morocco.
$$ E-mail\ address:\ alaouikarima2012@hotmail.fr$$}

\author{Najib Mahdou}
\address{Najib Mahdou\\Department of Mathematics, Faculty of Science and Technology of Fez, Box 2202,
University S.M. Ben Abdellah Fez, Morocco.
$$E-mail\ address:\ mahdou@hotmail.com$$}



\begin{abstract}
Let $f: A\rightarrow B$ be a ring homomorphism and let $J$ be an ideal of $B$. In this paper, we investigate the transfert of
the property of coherence to the amalgamation $A\bowtie^{f}J$. We provide necessary and sufficient conditions for  $A\bowtie^{f}J$ to be a coherent ring.

\end{abstract}

\maketitle

 \begin{section} {Introduction}
Throughout this paper, all rings are commutative with identity element, and all modules are unitary.\par
Let $R$ be a commutative ring. We say that an ideal is regular if
it contains a regular element, i.e; a non-zerodivisor element.\par
For a nonnegative integer $n$, an $R$-module $E$ is called $n$-presented
if there is an exact sequence of $R$-modules: \\
$$F_n \rightarrow  F_{n-1} \rightarrow \ldots F_1  \rightarrow F_0  \rightarrow E \rightarrow 0$$
 where each $F_i$ is a finitely generated free $R$-module. In particular, $0$-presented and $1$-presented $R$-modules are,
respectively, finitely generated and finitely presented
$R$-modules.\\

\par A ring $R$ is coherent if every finitely generated ideal of $R$
is finitely presented; equivalently, if $(0:a)$ and $I\cap J$ are
finitely generated for every $a\in R$ and any two finitely
generated ideals $I$ and $J$ of $R$. Examples of coherent rings
are Noetherian rings, Boolean algebras, von Neumann regular rings,
and Pr\"ufer/semihereditary rings. For instance see \cite{Gz2}. \\

\par Let $A$ and $B$ be two rings, let $J$ be an ideal of $B$ and let $f: A\rightarrow B$ be a ring homomorphism. In this setting, we can consider the
following subring of $A \times B$:\\
\begin{center}
$A\bowtie^{f}J = \{(a,f(a)+j) \diagup a \in A, j \in J\}$
\end{center}
called the amalgamation of $A$ with $B$ along $J$ with respect to
$f$ (introduced and studied by D'Anna, Finocchiaro, and Fontana in
\cite{AFF1,AFF2}). This construction is a generalization of the
amalgamated duplication of a ring along an ideal (introduced and
studied by D'Anna and Fontana in \cite{A, AF, AF2} and denoted by
$A \bowtie I$). Moreover, other classical constructions (such as
the $A+XB[X]$, $A+XB[[X]]$, and the $D+M$ constructions) can be
studied as particular cases of the amalgamation \cite[Examples 2.5
\& 2.6]{AFF1} and other classical constructions, such as the
Nagata's idealization and the CPI extensions (in the sense of
Boisen and Sheldon \cite{Boi}) are strictly related to it (see
\cite[Example 2.7 \& Remark 2.8]{AFF1}).\\

\par One of the key tools for studying $A\bowtie^{f}J$ is based on
the fact that the amalgamation can be studied in the frame of
pullback constructions \cite[Section 4]{AFF1}. This point of view
allows the authors in \cite{AFF1, AFF2} to provide an ample
description of various properties of $A\bowtie^{f}J$, in
connection with the properties of $A$, $J$ and $f$. Namely, in
\cite{AFF1}, the authors studied the basic properties of this
construction (e.g., characterizations for $A\bowtie^{f}J$  to be a
Noetherian ring, an integral domain, a reduced ring) and they
characterized those distinguished pullbacks that can be expressed
as an amalgamation.\\

\par This paper investigates a property of coherence in amalgamated
algebra along an ideal. Our results generate original examples
which enrich the current literature with new families of
non-Noetherian coherent rings.

 \end{section}

 \begin{section}{Main Results}
\bigskip


This section characterizes the amalgamated algebra along an ideal
$A\bowtie^{f}J$ to be a coherent ring. The main result (Theorem
2.2) examines the property of coherence that the amalgamation
$A\bowtie^{f}J$ might inherit from the ring A for some classes of
ideals J and homomorphisms $f$, and hence generates new examples
of non-Noetherian coherent rings.\\

 Let $f: A\rightarrow B$ be a ring homomorphism, $J$ be an
ideal of $B$ and let $n$ be a positive integer. Consider the
function $f^{n}: A^{n}\rightarrow B^{n}$ defined by
 $f^{n}((\alpha_{i})_{i=1}^{i=n}) = (f(\alpha_{i}))_{i=1}^{i=n}$. Obviously, $f^{n}$ is a ring homomorphism and $J^{n}$ is an ideal of $B^{n}$.
 This allows us to define $A^{n}\bowtie^{f^{n}}J^{n}$.\par
 Moreover, let $\phi: (A\bowtie^{f}J)^{n} \rightarrow
A^{n}\bowtie^{f^{n}}J^{n}$ defined by
$\phi((a_{i},f(a_{i})+j_{i})_{i=1}^{i=n}) =
((a_{i})_{i=1}^{i=n},f^{n}((a_{i})_{i=1}^{i=n})+(j_{i})_{i=1}^{i=n})$.
It is easily checked that $\phi$ is a ring isomorphism. So
$(A\bowtie^{f}J)^{n}$ and  $A^{n}\bowtie^{f^{n}}J^{n}$ are
isomorphic as rings.\par
Let $U$ be a submodule of $A^{n}$. Then $U \bowtie^{f^{n}}J^{n}:= \{(u,f^{n}(u)+j)\in A^{n} \bowtie^{f^{n}}J^{n} \diagup u \in U, j \in J^{n}\}$ is a
submodule of $A^{n} \bowtie^{f^{n}}J^{n}$. \\

\bigskip

Next, before we announce the main result of this section (Theorem
\ref{thm}), we make the following useful remark. \\

\bigskip

\begin{rem}
Let $f: A\rightarrow B$ be a ring homomorphism and let $J$ be an ideal of $B$. Then $f^{n}(\alpha a)= f(\alpha)f^{n}(a)$ for all
$\alpha \in A$ and $a \in A^{n}$.\\
Indeed, let $a= (a_{i})_{i=1}^{i=n} \in A^{n}$:\\
\begin{eqnarray}
\nonumber f^{n}(\alpha a ) &=& f^{n}(\alpha (a_{i})_{i=1}^{i=n}) =f^{n}((\alpha a_{i})_{i=1}^{i=n})\\
\nonumber &=&  (f(\alpha a_{i}))_{i=1}^{i=n}\\
\nonumber &=&(f(\alpha) f(a_{i}))_{i=1}^{i=n}\\
\nonumber &=& f(\alpha) (f(a_{i}))_{i=1}^{i=n}\\
\nonumber &=&f(\alpha) f^{n}((a_{i})_{i=1}^{i=n})= f(\alpha) f^{n}(a)
\end{eqnarray}
\end{rem}

\bigskip

 Now, to the main result:

\bigskip

\begin{thm}\label{thm}
Let $f: A\rightarrow B$ be a ring homomorphism and let $J$ be a proper ideal of $B$.
\begin{enumerate}
 \item  If $A\bowtie^{f}J$ is a coherent ring, then so is $A$.
 \item  Assume that $J$ and $f^{-1}(J)$ are finitely generated ideals of $f(A) + J$ and $A$ respectively. Then $A\bowtie^{f}J $ is a coherent ring if and only
 if $A$ and $f(A) + J$ are coherent rings.
 \item   Assume that  $J$ is a regular finitely generated ideal of $f(A) + J$. Then $A\bowtie^{f}J $ is a coherent ring if and only if $A$ and
 $f(A) + J$ are coherent rings and $f^{-1}(J)$ is a finitely generated ideal of $A$.
\end{enumerate}
\end{thm}


  \bigskip

Before proving Theorem \ref{thm}, we establish the following lemmas.

\bigskip

\begin{lem} \label{lem1}
Let $f: A\rightarrow B$ be a ring homomorphism and let $J$ be a
proper ideal of $B$. Then:
\begin{enumerate}
\item ${\{0\} \times J}$ (resp., $f^{-1}\{J\}\times \{0\}$) is a finitely generated ideal of $ A\bowtie^{f}J $ if and only if $J$
(resp., $f^{-1}\{J\}$) is a  finitely generated ideal of $f(A)+J$ (resp., $A$).
\item If $A\bowtie^{f}J$ is a coherent ring and $f^{-1}(J)$ is a finitely generated ideal of $A$, then $f(A) + J$ is a coherent ring.
\end{enumerate}

\end{lem}

\proof {\bf (1)} Assume that $J:= \sum_{i=1}^{i=n}(f(A)+J)k_{i}$
is a finitely generated ideal of $f(A) + J$, where $k_{i} \in J$.
It is clear that $\sum_{i=1}^{i=n} (A\bowtie^{f}J)
(0,k_{i})\subset {\{0\} \times J}$. Let $
x:=(0,\sum_{i=1}^{i=n}(f(\alpha_{i})+j_{i})k_{i}) \in {\{0\}
\times J}$, where $\alpha_{i} \in A$ and $j_{i} \in J$. Hence, $ x
= \sum_{i=1}^{i=n}(0,(f(\alpha_{i})+j_{i})k_{i})=
\sum_{i=1}^{i=n}(\alpha_{i}, f(\alpha_{i})+j_{i})(0,k_{i}) \in
\sum_{i=1}^{i=n}(A\bowtie^{f}J)((0,k_{i})$. Therefore, ${\{0\}
\times J} \subset \sum_{i=1}^{i=n} (A\bowtie^{f}J) (0,k_{i})$ and
so ${\{0\} \times J} = \sum_{i=1}^{i=n} (A\bowtie^{f}J)
(0,k_{i})$. Conversely, Assume that ${\{0\} \times J} :=
\sum_{i=1}^{i=n} (A\bowtie^{f}J) (0,k_{i})$ is a finitely
generated ideal of $ A\bowtie^{f}J $, where $k_{i} \in J$. It is
readily seen that $J= \sum_{i=1}^{i=n} (f(A)+J)k_{i}$, as
desired.\par Assume that $f^{-1}\{J\} := \sum_{i=1}^{i=n} A k_{i}$
is a finitely generated ideal of $A$, where $k_{i} \in
f^{-1}\{J\}$. It is obvious that $\sum_{i=1}^{i=n}A\bowtie^{f}J
(k_{i},0) \subset f^{-1}\{J\}\times \{0\} $. Let
$x=:(\sum_{i=1}^{i=n}\alpha_{i}k_{i},0) \in f^{-1}\{J\}\times
\{0\} $, where $ \alpha_{i} \in A$. Then
$x=\sum_{i=1}^{i=n}(\alpha_{i}k_{i},0) =
\sum_{i=1}^{i=n}(\alpha_{i},f(\alpha_{i}))(k_{i},0) \in
\sum_{i=1}^{i=n}(A\bowtie^{f}J)(k_{i},0)$. Therefore,
$f^{-1}\{J\}\times \{0\} \subset \sum_{i=1}^{i=n}(A\bowtie^{f}J)
(k_{i},0)$ and so $f^{-1}\{J\}\times \{0\} =
\sum_{i=1}^{i=n}A\bowtie^{f}J(k_{i},0)$. Conversely, Assume that
$f^{-1}\{J\}\times \{0\} := \sum_{i=1}^{i=n}
(A\bowtie^{f}J)(a_{i},0)$ is a finitely generated ideal of $
A\bowtie^{f}J $, where $a_{i} \in f^{-1}\{J\}$. It is easy to
check that $f^{-1}\{J\} = \sum_{i=1}^{i=n} Aa_{i}$, as desired.
\par {\bf (2)} Assume that $A\bowtie^{f}J $ is a coherent ring and
$f^{-1}\{J\}\times \{0\}$ is a finitely generated ideal of
$A\bowtie^{f}J $. Then $f(A) + J \cong
\frac{A\bowtie^{f}J}{f^{-1}\{J\}\times \{0\}}$ is a coherent ring
by \cite[Theorem 2.4.1]{Gz2}, as desired. \qed

\bigskip

\begin{lem}\label{lem2}
Let $f: A\rightarrow B$ be a ring homomorphism, $J$ be an ideal of $B$, and let $U$ be a submodule of $A^{n}$. Then:
\begin{enumerate}
 \item  Assume that $U$ is a finitely generated $A$-module and $J$ is a finitely generated ideal of $f(A)+J$. Then $U \bowtie^{f^{n}}J^{n}$ is a finitely generated $(A\bowtie^{f}J )$-module.
 \item  Assume that $f^{n}(U) \subset J^{n}$. Then $U \bowtie^{f^{n}}J^{n}$ is a finitely generated $(A\bowtie^{f}J )$-module if and only if $U$
 is a finitely generated $A$-module and $J$ is a finitely generated ideal of $f(A)+J$.
\end{enumerate}
\end{lem}

\proof {\bf (1)} Assume that $U :=\sum_{i=1}^{i=n}Au_{i}$ is a
finitely generated $A$-module, where $u_{i}\in U$ for all $
i\in\{1,.....n\}$ and $J^{n}:= \sum_{i=1}^{i=n}(f(A)+J) e_{i}$ is
a finitely generated $(f(A)+J)$-module, where $e_{i} \in J^{n}$
for all $ i\in\{1,.....n\}$. We claim that $U \bowtie^{f^{n}}J^{n}
= \sum_{i=1}^{i=n} (A\bowtie^{f}J)(u_{i},f^{n}(u_{i})) +
\sum_{i=1}^{i=n} (A\bowtie^{f}J)(0,e_{i})$. Indeed,
$\sum_{i=1}^{i=n} (A\bowtie^{f}J)(u_{i},f^{n}(u_{i})) +
\sum_{i=1}^{i=n} (A\bowtie^{f}J)(0,e_{i}) \subset U
\bowtie^{f^{n}}J^{n}$ since $(u_{i},f^{n}(u_{i}))\in U
\bowtie^{f^{n}}J^{n}$ for all $ i\in\{1,.....n\}$ and
$(0,e_{i})\in U \bowtie^{f^{n}}J^{n}$  for all $
i\in\{1,.....n\}$. Conversely, let $(x,f^{n}(x)+k)\in U
\bowtie^{f^{n}}J^{n}$, where $x \in U$ and $k \in J^{n}$. Hence,
$x =\sum_{i=1}^{i=n}\alpha_{i}u_{i} \in U$, for some
$\alpha_{i}\in A$ ($ i\in\{1,.....n\}$) and
$k=\sum_{i=1}^{i=n}(f(\beta_{i}) + j_{i})e_{i} \in J^{n}$, for
some $\beta_{i}\in A$ and $j_{i} \in J$ ($ i\in\{1,.....n\}$). We
obtain
\begin{eqnarray}
 \nonumber (x,f^{n}(x)+k) &=& (\sum_{i=1}^{i=n}\alpha_{i}u_{i},\sum_{i=1}^{i=n}f(\alpha_{i})f^{n}(u_{i})) + (0,\sum_{i=1}^{i=n}(f(\beta_{i}) + j_{i}) e_{i})\\
 \nonumber &=&  \sum_{i=1}^{i=n}(\alpha_{i},f(\alpha_{i}))(u_{i},f^{n}(u_{i})) + \sum_{i=1}^{i=n}(0,f(\beta_{i}) + j_{i})(0,e_{i})\\
\nonumber &=&
\sum_{i=1}^{i=n}(\alpha_{i},f(\alpha_{i}))(u_{i},f^{n}(u_{i})) +
\sum_{i=1}^{i=n}(\beta_{i},f(\beta_{i}) + j_{i})(0,e_{i}).
\end{eqnarray}
Consequently, $(x,f^{n}(x)+k)\in \sum_{i=1}^{i=n}
(A\bowtie^{f}J)(u_{i},f^{n}(u_{i})) + \sum_{i=1}^{i=n}
(A\bowtie^{f}J)(0,e_{i})$ since $(\alpha_{i},f(\alpha_{i})) \in
(A\bowtie^{f}J)$ for all $ i\in\{1,.....n\}$ and
$(\beta_{i},f(\beta_{i}) + j_{i}) \in A\bowtie^{f}J$ for all $
i\in\{1,.....n\}$ and hence $U \bowtie^{f^{n}}J^{n} =
\sum_{i=1}^{i=n} (A\bowtie^{f}J)(u_{i},f^{n}(u_{i})) +
\sum_{i=1}^{i=n} (A\bowtie^{f}J)(0,e_{i})$ is a finitely generated
$(A\bowtie^{f}J )$-module, as desired.
\par {\bf (2)} Assume that $f^{n}(U) \subset J^{n}$. If $U$ is a
finitely generated $A$-module and $J$ is a finitely generated
ideal of $f(A)+J$, then $U \bowtie^{f^{n}}J^{n}$ is a finitely
generated $(A\bowtie^{f}J )$-module by (1). Conversely, assume
that $U \bowtie^{f^{n}}J^{n} := \sum_{i=1}^{i=n} (A\bowtie^{f}J
)(u_{i}, f^{n}(u_{i})+ k_{i})$ is a finitely generated
$(A\bowtie^{f}J )$-module, where, $u_{i} \in U$ and $k_{i} \in
J^{n}$ for all $ 1 \leq i \leq  n$. It is clear that $U =
\sum_{i=1}^{i=n} Au_{i}$. On the other hand, we claim that $J^{n}
= \sum_{i=1}^{i=n}(f(A)+J)(f^{n}(u_{i}) +k_{i})$. Indeed, let $j
\in J^{n}$. Then $(0,j) = \sum_{i=1}^{i=n}
(\alpha_{i},f(\alpha_{i}) + j_{i})(u_{i}, f^{n}(u_{i})+ k_{i})$
for some $\alpha_{i} \in A$ and $j_{i}\in J$. So
 $j = \sum_{i=1}^{i=n}(f(\alpha_{i}) + j_{i})(f^{n}(u_{i})+ k_{i}) \in \sum_{i=1}^{i=n}(f(A)+J)(f^{n}(u_{i}) +k_{i})$. Thus
 $J^{n} \subset \sum_{i=1}^{i=n}(f(A)+J)(f^{n}(u_{i}) +k_{i})$. But  $f^{n}(u_{i}) \in J^{n}$ for all $i = 1,...n$ since
 $f^{n}(U) \subset J^{n}$. Hence,
 $(f^{n}(u_{i}) + k_{i}) \in J^{n}$ $\forall i$ and so $\sum_{i=1}^{i=n}(f(A)+J)(f^{n}(u_{i}) +k_{i})\subset J^{n}$.
 Therefore, $J^{n}=\sum_{i=1}^{i=n}(f(A)+J)(f^{n}(u_{i}) +k_{i}) $ is a finitely generated $(f(A)+J)$-module and so $J$ is a finitely
 generated ideal of $(f(A)+J)$, completing the proof of Lemma \ref{lem2}.
\qed

\bigskip


\par At this point, it is worthwhile recalling that an $R$-module $M$ is called a coherent $R$-module if it is finitely generated and every finitely
generated submodule of $M$ is finitely presented.\\

\begin{lem}\label{lem3}
Let $f: A\rightarrow B$  be a ring homomorphism, and $J$ be an
ideal of $B$. Assume that $J$ and $f^{-1}(J)$ are finitely
generated ideals of  $f(A) + J$ and $A$ respectively. Then
$f^{-1}\{J\}\times \{0\}$ is a coherent  $(A\bowtie^{f}J )$-module
provided $A$ is a coherent ring.
\end{lem}

\proof Since $f^{-1}\{J\}\times \{0\}$ is a finitely generated
($A\bowtie^{f}J$)-module, it remains to show that every finitely
generated submodule of $f^{-1}\{J\}\times \{0\}$ is finitely
presented. Assume that $A$ is a coherent ring and let $N$ be a
finitely generated submodule of $f^{-1}\{J\}\times \{0\}$. It is
clear that $N = I \times \{0\}$, where $I = \sum_{i=1}^{i=n} A
a_{i}$ for some positive integer $n$ and  $a_{i} \in I$. Consider
the exact sequence of $A$-modules:
\begin{center}
$ 0\rightarrow Kerv\rightarrow A^{n} \rightarrow I \rightarrow 0 $ \hspace{2 cm} (1)
\end{center}
where $v((\alpha_{i})_{i=1}^{i=n}) = \sum_{i=1}^{i=n}\alpha_{i}a_{i}$. Then $Kerv = \{(\alpha_{i})_{i=1}^{i=n} \in A^{n} \diagup \sum_{i=1}^{i=n}\alpha_{i}a_{i} = 0 \}$.
On the other hand, it is easily verified that $N = \sum_{i=1}^{i=n} A\bowtie^{f}J(a_{i},0)$. Consider  the exact sequence of $(A\bowtie^{f}J)$-modules:\\
\begin{center}
$ 0\rightarrow Keru\rightarrow (A\bowtie^{f}J)^{n} \rightarrow N \rightarrow 0 $\hspace{2 cm} (2)
\end{center}
where $u((\alpha_{i},f(\alpha_{i}) + j_{i})_{i=1}^{i=n}) =
\sum_{i=1}^{i=n}(\alpha_{i},f(\alpha_{i}) + j_{i})(a_{i},0)$.
Then, $Keru = \{(\alpha_{i},f(\alpha_{i}) + j_{i})_{i=1}^{i=n} \in
(A\bowtie^{f}J)^{n} \diagup \sum_{i=1}^{i=n}\alpha_{i}a_{i} =
0\}$. So $Keru =$ \\
$\{((\alpha_{i})_{i=1}^{i=n},f^{n}((\alpha_{i})_{i=1}^{i=n}) +
(j_{i})_{i=1}^{i=n}) \in A^{n}\bowtie^{f^{n}}J^{n} \diagup
(\alpha_{i})_{i=1}^{i=n} \in Kerv\}$ and hence $Keru = Kerv
\bowtie^{f^{n}}J^{n}$. But $I$ is a finitely presented ideal of
$A$ since $A$ is a coherent ring, so $Kerv$ is a finitely
generated $A$-module (by a sequence (1)) and hence $Keru = Kerv
\bowtie^{f^{n}}J^{n}$ is a finitely generated
$(A\bowtie^{f}J)$-module (by lemma \ref{lem2} (1)). Therefore, N
is a finitely presented $(A\bowtie^{f}J)$-module by a sequence (2)
and hence $f^{-1}\{J\}\times \{0\}$ is a coherent $A\bowtie^{f}J
$-module, to complete the proof of Lemma \ref{lem3}. \qed

\bigskip

\begin{lem} \label{lem4}
Let $f: A\rightarrow B$ be a ring homomorphism, and $J$ be an ideal of $B$. If $A\bowtie^{f}J $ is a coherent ring and $J$ is a regular
ideal of $f(A) + J$, then $f^{-1}(J)$ is a finitely generated ideal of $A$.
\end{lem}

\proof
Assume that $A\bowtie^{f}J $ is a coherent ring and $J$ contains a regular element $k$. Set $c = (0,k) \in A\bowtie^{f}J$. One can easily check that:
\begin{eqnarray}
 \nonumber (0:c) &=& \{(a,f(a) + j) \in A\bowtie^{f}J \diagup (a,f(a) + j)(0,k)=0\} \\
 \nonumber  &=&  \{(a,f(a) + j) \in A\bowtie^{f}J \diagup (f(a) + j) k = 0 \}\\
\nonumber   &=& \{(a,f(a) + j) \in A\bowtie^{f}J \diagup f(a) + j = 0\}\\
\nonumber   &=&\{(a,0) \in A\bowtie^{f}J \diagup a \in f^{-1}\{J\}\}\\
\nonumber  &=& f^{-1}\{J\}\times \{0\}.
\end{eqnarray}
Since $A\bowtie^{f}J $ is a coherent ring, then $(0:c)=
f^{-1}\{J\}\times \{0\}$ is a finitely generated ideal of
$A\bowtie^{f}J $. Therefore, $f^{-1}\{J\}$ is a finitely generated
ideal of $A$, as desired. \qed
\bigskip
\bigskip

{\bf Proof of Theorem 2.2} \\
\par {\bf (1)} If $A\bowtie^{f}J $ is a coherent ring, then $A$ is
a coherent ring by \cite[Theorem 4.1.5]{Gz2} since $A$ is a module
retract of $A\bowtie^{f}J $.
\par {\bf (2)} Assume that $J$ and $f^{-1}(J)$ are finitely generated
ideals of  $f(A) + J$ and $A$ respectively. Then $A$ and $f(A) +
J$ are coherent rings since $A\bowtie^{f}J$ is a coherent ring (by
Theorem \ref{thm} (1) and Lemma \ref{lem1} (2)). Conversely,
assume that $A$ and $f(A)+J$ are coherent rings. Since
$\frac{A\bowtie^{f}J}{f^{-1}\{J\}\times \{0\}} \cong f(A) + J$,
$f(A) + J$ is a coherent ring and $f^{-1}\{J\}\times \{0\}$ is a
coherent $A\bowtie^{f}J $-module (by Lemma \ref{lem3}), then $
A\bowtie^{f}J $ is a coherent ring (by \cite[Theorem 2.4.1]{Gz2}).
\par {\bf (3)} Follows immediately from Theorem \ref{thm} (2)
and Lemma \ref{lem4}. This completes the proof of the main
Theorem. \qed
\bigskip
\bigskip

The following Corollary is an immediate consequence of Theorem
\ref{thm} (3).

\begin{cor} \label{cor1}
Let $f: A\rightarrow B$ be a ring homomorphism, $B$ be an integral
domain and let $J$ be a proper and finitely generated ideal of
$f(A)+J$. Then $A\bowtie^{f}J $ is a coherent ring if and only if
$A$ and $f(A) + J$ are coherent rings and $f^{-1}(J)$ is a
finitely generated ideal of $A$.
\end{cor}
\bigskip

The Corollary below follows immediately from Theorem \ref{thm} (2)
which examines the case of the amalgamated duplication.

\begin{cor} \label{cor2}
Let $A$ be a ring and $I$ be a proper ideal of A.
\begin{enumerate}
\item If $A\bowtie I$ is a coherent ring, then so is A.
\item Assume that $I$ is a finitely generated ideal of $A$. Then $A\bowtie I$ is a coherent ring if and only if $A$  is a coherent ring.
   \end{enumerate}

\end{cor}

\bigskip


The next Corollary is an immediate consequence of Theorem
\ref{thm} (2).

\begin{cor} \label{cor3}
Let $A$ be a ring, $I$ be an ideal of $A$, $B:=\frac{A}{I}$, and let $f: A\rightarrow B$ be the canonical homomorphism ($f(x) = \overline{x}$).
\begin{enumerate}
\item Assume that $J$ and $f^{-1}(J)$ are finitely generated ideals of $B$ and $A$ respectively. Then $A\bowtie^{f}J $ is a coherent ring if and only if $A$ and $B$ are coherent rings.
 \item   Assume that  $J$ is a regular finitely generated ideal of $B$. Then $A\bowtie^{f}J $ is a coherent ring if and only if $A$ and $B$ are coherent rings and $f^{-1}(J)$ is a finitely generated ideal of $A$.
\end{enumerate}
\end{cor}

\bigskip

 The aforementioned result enriches the literature with new
examples of coherent rings which are non-Noetherian rings.

\begin{exam}
Let $A$ be a non-Noetherian coherent ring, $I$ be a finitely
generated ideal of $A$, $f: A\rightarrow B (= \frac{A}{I})$ be the
canonical homomorphism, and let $J$ be a finitely generated ideal
of $A$. Then $A\bowtie^{f}\overline{J} $ is a non-Noetherian
coherent ring.
\end{exam}

\proof By Corollary \ref{cor3}, $A\bowtie^{f}\overline{J} $
is a coherent ring since $A$ and $B$ are both coherent rings and
$J$ is a finitely generated ideal of $A$. On the other hand,
  $A\bowtie^{f}\overline{J}
$ is a non-Noetherian ring by \cite[Proposition 5.6, p. 167]{AFF1}
since $A$ is a non-Noetherian ring. \qed

\bigskip

\begin{exam}
Let $K$ be a field and consider the power series ring $A = K[[X_{1},.....,X_{n},...]]$ and let $I:= X_{1}A+X_{2}A$. Then
$A\bowtie I$ is a non-Noetherian coherent ring.\\
\end{exam}

\proof By Corollary \ref{cor2}, $A\bowtie I $ is a coherent ring
since $A$ is coherent and $I$ is a finitely generated ideal of
$A$. On the other hand, $A\bowtie I $ is a non-Noetherian ring by
\cite[Corollary 3.3]{AF2} since $A$ is a non-Noetherian ring. \qed

\bigskip

\begin{exam}
Let $A:=\Z+X\Q[X]$, where $\Z$ is the ring of integers, and $\Q$
is the field of rational numbers. Let $I:=X\Q[X]$,
$B:=\frac{A}{I}(\cong \Z) $, $f: A\rightarrow B $ be the canonical
homomorphism and let $J$ be a nonzero ideal of $B$. Then
$A\bowtie^{f}J $ is a non-Noetherian coherent ring.
\end{exam}

\proof By Corollary \ref{cor3}, $A\bowtie^{f}J$ is a coherent ring
since $A$ and $B$ are both coherent rings and $J$ (resp.,
$f^{-1}(J) = n_{0}\Z +X\Q[X] $ for some positive integer $n_{0}$)
is a finitely generated ideal of $B$ (resp., $A$). On the other
hand, $A\bowtie^{f}J$ is a non-Noetherian ring by
\cite[Proposition 5.6, p. 167]{AFF1} since $A$ is a non-Noetherian
ring. \qed
\end{section}

\bigskip
\bigskip




\bigskip\bigskip

\end{document}